\documentclass[10pt]{article}

\usepackage[latin1]{inputenc}
\usepackage{url, amsmath,enumerate,fancyhdr,amssymb, amsthm, 
pstricks, epsf, lscape, ifpdf, graphicx, float}
\usepackage[margin=2.5cm]{geometry}
\usepackage{mathtools}

\mathchardef\mhyphen="2D

\mathchardef\mhyphen="2D

\newcommand\clps{\hspace{.9cm}}
\newcommand\clpss{\hspace{0.4 cm}}
\newcommand\paux{(P_{\rm aux})}
\newcommand\daux{(D_{\rm aux})}
\newcommand\phom{(P_{\rm hom})}
\newcommand\dhom{(D_{\rm hom})}  
\newcommand\pgen{(P_{\rm gen})}
\newcommand\dgen{(D_{\rm gen})}  

\newtheorem{Definition}{Definition}
\newtheorem{Proposition}{Proposition}
\newtheorem{Lemma}{Lemma}
\newtheorem{Theorem}{Theorem}
\newtheorem{Corollary}{Corollary}

\newtheorem{Assumption}{Assumption}

\newtheorem{Remark}{Remark}
\newtheorem{Example}{Example}

\newcommand{\INFEAS}{\operatorname{INFEAS}}
\newcommand{\FEAS}{\operatorname{FEAS}}

\newcommand{\ri}{\operatorname{ri}}

\newcommand{\eps}{\epsilon} 
\newcommand{\bpx}{\begin{pmatrix}}
\newcommand{\epx}{\end{pmatrix}}
\newcommand{\bpxastr}{\begin{pmatrix*}[r]}
\newcommand{\epxastr}{\end{pmatrix*}}
\newcommand{\bbx}{\begin{bmatrix}}
\newcommand{\ebx}{\end{bmatrix}}

\newcommand{\bdef}{\begin{Definition}} 
\newcommand{\commentout}[1]{}
\newcommand{\co}[1]{}

\newcommand{\nin}{\noindent}

\newcommand{\pf}[1]{\vspace{.35cm} \nin {\bf Proof {#1} }}

\newcommand{\sym}[1]{{\cal S}^{#1}}
\newcommand{\psd}[1]{{\cal S}_+^{#1}}
\newcommand{\pd}[1]{{\cal S}_{++}^{#1}}
\newcommand{\rad}[1]{\mathbb{R}^{#1}}

\newcommand{\eref}[1]{(\ref{#1})}

\newcommand{\A}{ {\cal A} }
\newcommand{\B}{ {\cal B} }

\newcommand{\la}{\langle}
\newcommand{\ra}{\rangle}

\newcommand{\beq}{\begin{equation}}
\newcommand{\eeq}{\end{equation}}
\newcommand{\beqa}{\begin{eqnarray}}
\newcommand{\eeqa}{\end{eqnarray}}
\newcommand{\ba}{\begin{array}}
\newcommand{\ena}{\end{array}}
\newcommand{\bac}{\begin{array}{ccccccccccc}}
\newcommand{\eac}{\end{array}}
\newcommand{\bprop}{\begin{Proposition}}
\newcommand{\eprop}{\end{Proposition}}

\newcommand{\conp}{{\rm co \mhyphen}{\cal NP}}

\newcommand{\beqast}{\begin{eqnarray*}}
\newcommand{\eeqast}{\end{eqnarray*}}
\newcommand{\benum}{\begin{enumerate}}
\newcommand{\eenum}{\end{enumerate}}
\newcommand{\bit}{\begin{itemize}}
\newcommand{\eit}{\end{itemize}}
\newcommand{\bth}{\begin{Theorem}}
\newcommand{\enth}{\end{Theorem}}
\newcommand{\ble}{\begin{Lemma}}
\newcommand{\ele}{\end{Lemma}}
\newcommand{\bex}{\begin{Example}}
\newcommand{\eex}{\end{Example}}
\newcommand{\bcor}{\begin{Corollary}}
\newcommand{\ecor}{\end{Corollary}}
\newcommand{\brem}{\begin{Remark}}
\newcommand{\erem}{\end{Remark}}
\newcommand{\bass}{\begin{Assumption}}
\newcommand{\eass}{\end{Assumption}}

\newcommand{\LRA}{\Leftrightarrow}

\newcommand{\bsmx}{\begin{small} \begin{pmatrix}}
\newcommand{\esmx}{\end{pmatrix} \end{small}}

\title{\bf Exact duality in semidefinite programming based on elementary reformulations \footnotemark[1]}

\author{Minghui Liu \hspace{2.5cm} G\'{a}bor Pataki  \\ {\bf minghui@unc.edu \hspace{1.6cm}\bf gabor@unc.edu} \vspace{.5cm} \\ Department of Statistics and Operations Research \\
University of North Carolina at Chapel Hill }

\begin{document}
\maketitle

\renewcommand{\thefootnote}{\fnsymbol{footnote}} \footnotetext[1]{The paper's previous title was "A short proof of infusibility and generating all infeasible semidefinite programs"}

\begin{abstract}
In semidefinite programming (SDP), 
unlike in linear programming, Farkas' lemma may fail to prove infeasibility. Here 
we obtain an exact, short certificate of infeasibility in SDP by an elementary approach:
we reformulate any semidefinite system of the form 
\beq\label{p} \tag{P}
\ba{rcl}
A_{i}\bullet X & = & b_{i} \,\,(i=1,\ldots,m) \\
X & \succeq & 0.
\ena
\eeq
using only elementary row operations, and rotations. 
When $(P)$ is infeasible, the reformulated system is trivially infeasible.  
When $(P)$ is feasible, the reformulated system has strong duality with its Lagrange dual 
for all objective functions. 
As a corollary, we obtain algorithms to generate the constraints of {\em all} infeasible SDPs 
and the constraints of {\em all} feasible SDPs with a fixed rank maximal solution. 

We give two methods to construct our elementary reformulations.
One is direct, and based on a simplified facial reduction algorithm, 
and the other is obtained by adapting the facial reduction algorithm of 
Waki and Muramatsu. 

In somewhat different language, our reformulations provide 
a standard form of spectrahedra, 
to easily verify either their emptiness, or 
a tight upper bound on the rank of feasible solutions. 
\end{abstract}



{\em OR/MS subject classification:} Primary: convexity; secondary: programming-nonlinear-theory

\pagestyle{myheadings}
\thispagestyle{plain}

\section{Introduction. The certificate of infeasibility and its proof}

Semidefinite programs (SDPs) naturally generalize linear programs and 
share some of the  duality theory of linear programming. 
However, the  value of an SDP may not be attained, it may differ from the 
value of its Lagrange dual, and the simplest version of 
Farkas' lemma may fail to prove infeasibility in semidefinite programming. 

Several alternatives of the traditional Lagrange dual, and Farkas' lemma are known, which we will review in detail 
below: see 
Borwein and Wolkowicz \cite{BorWolk:81, BorWolk:81B}; 
Ramana \cite{Ramana:97}; 
Ramana, Tun\c{c}el, and Wolkowicz \cite{RaTuWo:97}; 
Klep and Schweighofer \cite{KlepSchw:12}; Waki and Muramatsu 
\cite{WakiMura:12}, and the second author \cite{Pataki:13}. 

We consider semidefinite systems  of the form \eref{p}, 
where the $A_i$ are $n$ by $n$ 
symmetric matrices, the $b_i$ scalars, $X \succeq 0$ means that $X$ is symmetric, 
positive semidefinite (psd), 
and the $\bullet$ dot product of symmetric matrices is 
the trace of their regular product. To motivate our results on infeasibility, 
we consider the 
instance 
\beq \label{mot-ex}
\ba{cclcc}
\bpx 1 & 0 & 0 \\
     0 & 0 & 0 \\
     0 & 0 & 0 
\epx &\bullet&  X & = & 0 \\
\vspace{0.2cm}
\bpx 0 & 0 & 1 \\
     0 & 1 & 0 \\
     1 & 0 & 0 
\epx \, &\bullet& \, X & = & -1 \\
                  && X & \succeq & 0,
\ena
\eeq
which is trivially infeasible: to see why, suppose that $X = (x_{ij})_{i,j=1}^3$ is feasible in it.
Then $x_{11}=0, $ hence the first row and column of $X$ are zero
by psdness, so  the second constraint implies $x_{22} = -1, \,$ which is a contradiction.
Thus the internal structure of the system itself proves its infeasibility.

The goal of this short note is twofold. In Theorem \ref{ref-SDP} we show that 
a basic transformation reveals such a simple structure -- which proves infeasibility -- 
in {\em every} infeasible semidefinite system. 
For feasible systems we give a similar reformulation -- in Theorem \ref{ref-SDP-2} --
which trivially has strong duality with its Lagrange dual for all objective functions.
\bdef
We obtain an {\em elementary semidefinite (ESD-) reformulation, or elementary reformulation} 
of \eref{p} by applying a sequence of the following operations:
\bit
\item[(1)] Replace $(A_j, b_j)$ by $(\sum_{i=1}^{m}y_{i}A_{i}, \sum_{i=1}^{m}y_{i}b_{i})$, 
where $y \in \rad{m}, \, y_j \neq 0.$ 
\item[(2)] Exchange two equations. 
\item[(3)] Replace $A_i$ by $V^T A_i V$ for all $i, \,$ where 
$V$ is an invertible matrix.
\eit
\end{Definition}
ESD-reformulations clearly preserve feasibility.
Note that operations (1) and (2) are also used in Gaussian elimination: 
we call them elementary row operations (eros). 
We call operation (3) a rotation. 
Clearly, we can assume that a rotation is applied only once, when 
reformulating \eref{p}; then $X$ is feasible for \eref{p} if and only if
$V^{-1} X V^{-T}$ is feasible for the reformulation.

\bth \label{ref-SDP}
The system \eref{p} is infeasible, if and only if it has an elementary semidefinite reformulation of the form
\begin{equation}\label{pref}\tag{P$_{\rm ref}$}
\begin{array}{rcl} 
 A_{i}^{\prime}\bullet  X & = & 0 \, (i=1,\ldots,k) \\
 A_{k+1}^{\prime}\bullet X & = & -1\\
 A_{i}^{\prime}\bullet X & = & b_{i}^{\prime} \, (i=k+2,\ldots,m) \\
 X & \succeq & 0
\end{array}
\end{equation}
where $k \geq 0, \,$ and the $A_i^\prime$ are of the form 

$$
A_{i}^{\prime} = 
\bordermatrix{
& \overbrace{\qquad \qquad \qquad}^{\textstyle r_{1}+\ldots+r_{i-1}} & \overbrace{\qquad}^{\textstyle r_{i}} & \overbrace{\qquad \qquad \qquad\quad}^{\textstyle n-r_{1}-\ldots-r_{i}} \cr\\
& \times  &  \times  &  \times \cr
& \times  &  I   &  0 \cr
& \times  &  0  &  0 \cr} 
$$
\noindent for $i=1,\ldots,k+1$, with $r_{1},\ldots,r_{k}>0, r_{k+1}\geq 0, \,$ the $\times$ symbols 
correspond to blocks with arbitrary elements, and matrices $A_{k+2}^\prime, \dots, A_{m}^\prime$ 
and scalars $b_{k+2}^\prime, \dots, b_{m}^\prime$ are arbitrary.  
\enth

To motivate the reader, we now give a very simple, full proof of the ``if'' direction.
It suffices to prove that \eref{pref} is infeasible, so 
assume to the contrary that $X$ is 
feasible in it. The constraint 
$A_{1}^{\prime}\bullet X=0$ and $X\succeq 0$ implies that the upper left 
$r_1$ by $r_1$ block of $X$ is  zero, and $X \succeq 0$ 
proves that the first $r_1$ rows and columns 
of $X$ are zero. Inductively, from the first $k$ constraints  we deduce 
that the first $\sum_{i=1}^{k}r_{i}$ rows and columns of $X$ are zero.

Deleting the first $\sum_{i=1}^{k}r_{i}$ rows and columns from $A_{k+1}^{\prime}$ we obtain 
a psd matrix, hence 
$$
\ba{rclcl} 
A_{k+1}^{\prime}\bullet X & \geq & 0, 
\end{array}
$$
contradicting the $(k+1)^{st}$ constraint in \eref{pref}.
\qed

Note that Theorem \ref{ref-SDP} allows us to systematically generate {\em all}
infeasible semidefinite systems: to do so, we only need to generate systems of the
form \eref{pref}, and reformulate them.
We comment more on this in Section \ref{conclusion}. 

We now review relevant literature in detail, and its connection to our results.
For surveys and textbooks on SDP, we refer to  Todd 
\cite{Todd:00}; Ben-Tal and Nemirovskii \cite{BentalNem:01}; Saigal et al \cite{SaigVandWolk:00};
Boyd and Vandenberghe \cite{BoydVand:04}. For 
treatments of their duality theory see 
Bonnans and Shapiro \cite{BonnShap:00}; Renegar \cite{Ren:01} and G$\ddot{\mathrm{u}}$ler \cite{Guler:10}. 

The fundamental facial reduction algorithm of 
Borwein and Wolkowicz \cite{BorWolk:81, BorWolk:81B} ensures strong 
duality in a  possibly nonlinear conic system by replacing the underlying cone by a suitable face.
Ramana in \cite{Ramana:97} constructed an extended strong dual for SDPs, 
which uses $O(n)$ copies of the original system, and extra variables. His dual 
leads to an exact Farkas' lemma. Though these approaches seem at first quite different,
Ramana, Tun\c{c}el, and Wolkowicz in \cite{RaTuWo:97} proved the correctness of 
Ramana's dual from the algorithm in \cite{BorWolk:81, BorWolk:81B}.

The algorithms in \cite{BorWolk:81, BorWolk:81B} assume that the system is feasible.
The simplified algorithm of Waki and Muramatsu  in \cite{WakiMura:12}, which works for conic linear systems,
disposes with this assumption, and allows one to prove infeasibility.
We state here that our reformulations can be obtained by suitably modifying the 
algorithm in \cite{WakiMura:12}; 
we describe the connection in detail in Section \ref{conclusion}. At the same time we provide a direct, and 
entirely elementary construction. 

More recently, Klep and Schweighofer in \cite{KlepSchw:12} 
proposed a strong dual and exact Farkas' lemma 
for SDPs. Their dual resembles Ramana's; 
however, it is based on ideas from algebraic geometry, namely sums of squares 
representations, not convex analysis. 

The second author in \cite{Pataki:13} described a simplified facial reduction algorithm,  
and generalized Ramana's dual to conic linear systems over {\em nice} cones
(for literature on nice cones, see \cite{ChuaTuncel:08}, \cite{Vera:13}, \cite{Pataki:12}).
We refer to P\'{o}lik and Terlaky \cite{PolikTerlaky:09}
for a generalization of Ramana's dual for conic LPs 
over homogeneous cones. 
Elementary reformulations of semidefinite systems first appear in 
\cite{Pataki:10}. There the second author uses them to bring a system into a form to easily check 
whether it has strong duality with its dual for all objective functions.

Several papers -- see for instance P\'{o}lik and Terlaky 
\cite{PolikTerlaky:09b} on stopping criteria for conic optimization --
point to the need of having more infeasible instances and we hope that our results 
will be useful in this respect. 
In more recent related work, 
Alfakih \cite{Alfakih:14} gave a certificate of 
the maximum rank in a feasible semidefinite system, using
a sequence of matrices, somewhat similar to the constructions in the duals of 
\cite{Ramana:97, KlepSchw:12}, and used it 
in an SDP based proof of a result of Connelly and Gortler on rigidity 
\cite{Connelly:14}. Our Theorem \ref{ref-SDP-2} gives 
such a certificate using elementary reformulations. 

We say that an infeasible SDP is weakly infeasible, if the 
traditional version of Farkas' lemma fails to prove its infeasibility.
We refer to Waki \cite{Waki:12} for a systematic method to generate weakly infeasible SDPs from Lasserre's 
relaxation of polynomial optimization problems; and to Lourenco et al. \cite{Lourenco:13} 
for an error-bound based reduction procedure to simplify weakly infeasible SDPs.

We organize the rest of the paper as follows. After introducing notation, 
we describe an algorithm to 
find the reformulation \eref{pref}, and a  constructive proof of 
the ``only if'' part of Theorem \ref{ref-SDP}. 
The algorithm is based on facial reduction; however, it is simplified so we do not
need to explicitly refer to faces of the semidefinite cone. 
The algorithm needs a subroutine to solve a 
primal-dual pair of SDPs. 
In the SDP pair the primal will always be strictly feasible, but the 
dual possibly not,  and we need to solve them in exact arithmetic.
Hence our algorithm may not run in polynomial time. At the same time it is quite simple, 
and we believe that it will be useful to verify the infeasibility of small instances.
We then illustrate the algorithm with Example \ref{ex1}. 

In Section \ref{feasible} we 
present our reformulation of feasible systems. 
Here we modify our algorithm 
to construct the reformulation \eref{pref} (and hence detect infeasibility); 
or to construct a reformulation that 
is easily seen to have strong duality with its Lagrange dual for all objective 
functions. 

We denote by $\sym{n}, \, \psd{n}, \,$ and $\pd{n}$ the set of 
symmetric, symmetric psd, and symmetric positive definite (pd) matrices of order $n, \,$ respectively. 
For a closed, convex cone $K$ we write $x \geq_K y$ to denote 
$x - y \in K, \,$ and denote the relative interior of $K$ by $\ri K, \,$ and its 
dual cone by $K^*, \,$ i.e.,
$$
\ba{rcl}
K^* & = & \{ \, y \, | \, \la x, y \ra \geq 0 \, \forall x \in K \, \}.
\ena
$$
For some $p < n$ we denote by $0 \oplus \psd{p}$ 
the set of $n$ by $n$ matrices
with the lower right $p$ by $p$ 
corner psd, and the rest of the components zero.
If $K = 0 \oplus \psd{p}, \,$ then $\ri K = 0 \oplus \pd{p}, \,$ and 
$$
K^* \, = \, \biggl\{ \, \bpx Z_{11} & Z_{12} \\
                                     Z_{12}^T & Z_{22} \epx \, : \, Z_{22} \in \psd{p} \biggr\}.
$$
For a matrix $Z \in K^*$ partitioned as above, and $Q \in \rad{p \times p}$ 
we will use the formula 
\beq \label{rotatez} 
\bpx I_{n-p} & 0 \\ 0 & Q \epx^T Z \bpx I_{n-p} & 0 \\ 0 & Q \epx \, = \, \bpx Z_{11} & Z_{12} Q \\ Q^{T} Z_{12}^T & Q^T Z_{22} Q \epx 
\eeq
in the reduction step of our algorithm that converts \eref{p} into 
\eref{pref}: we will choose 
$Q$ to be full rank, so that $Q^T Z_{22} Q$ is diagonal. 

We will rely on the following general conic linear system: 
\begin{equation}\label{pgen}
\begin{array}{rcl} 
\A(x) & = & b \\
\B(x) & \leq_K & d, 
\end{array}
\end{equation}
where $K$ is a closed, convex cone, and $\A$ and $\B$ are linear operators, 
and consider the primal-dual pair of conic LPs
$$
\ba{rrlcrrcll}
&\sup & \la c, x\ra & \clps & \inf & \la b, y\ra+\la d, z\ra & & &\\
\pgen & s.t. & x \text{ is feasible in (1.3)} & \clps & s.t. & \A^{*}(y)+\B^{*}(z) & = & c & \dgen \\
& & & \clps & & z & \in & K^{*}&
\ena
$$
where $\A^*$ and $\B^*$ are the adjoints of $\A$ and $\B, \,$ respectively. 

\begin{Definition}
We say that 
\benum
\item strong duality holds between $\pgen$ and $\dgen$, 
if their optimal values agree, and the latter value is attained, when finite;
\item \eref{pgen} is well behaved, if strong duality holds betweeen $\pgen$ and $\dgen$ for {\em all} $c$ objective functions;
\item \eref{pgen} is strictly feasible, if $d - \B(x) \in \ri K \,$ for some feasible $x.$
\eenum
\end{Definition}
We will use the following lemma: 

\ble \label{lemma-strong}
If \eref{pgen} is strictly feasible, or $K$ is polyhedral, then \eref{pgen} is well behaved.
\ele

When $K = 0 \oplus \psd{p}$ for some $p \geq 0, \,$ then $\pgen\mhyphen\dgen$ are a primal-dual pair of SDPs.
To solve them efficiently, we must assume that both are 
strictly feasible; strict feasibility of the latter means
that there is a feasible $(y, z)$ with  $z \in \ri K^*.$ 

The system \eref{p} is trivially infeasible, if the 
alternative system below is feasible:
\beq\label{palt} 
\ba{rcl}  
y & \in & \rad{m} \\
\sum_{i=1}^{m}y_{i}A_{i} & \succeq & 0\\
\sum_{i=1}^{m}y_{i} b_{i} & =       & -1;
\ena
\eeq
in this case we say that \eref{p} is strongly infeasible.
Note that system \eref{palt} generalizes Farkas' lemma from linear programming.
However, \eref{p} and \eref{palt} may both be infeasible, 
in which case we say that  \eref{p} is weakly infeasible. For instance, 
the system \eref{mot-ex} is weakly infeasible.

\pf{of ''only if'' in Theorem \ref{ref-SDP}}
The proof relies only on Lemma \ref{lemma-strong}.
We start with the system \eref{p}, which we assume to be infeasible.

In a general step we have a system
\beq \label{pprime}\tag{P$^\prime$}
\ba{rcl}
A_{i}^{\prime} \bullet X & = & b_{i}^{\prime} \, (i=1,\ldots,m) \\
X & \succeq & 0, 
\ena
\eeq
where for some $\ell \geq 0$ and $r_1 > 0, \dots, r_{\ell} > 0$ 
the $A_{i}^{\prime}$ matrices are as  required by Theorem 
\ref{ref-SDP}, and $b_{1}^{\prime} \, = \, \dots \, = \, b_{\ell}^{\prime} = 0. \,$ 
At the start 
$\ell = 0, \,$ and in a general step we have $0 \leq \ell < \min \{n, m \}.$ 

Let us define 
$$
r := r_1 + \dots + r_{\ell}, \, K := 0 \oplus \psd{n-r}, 
$$
and note that if $X \succeq 0$ satisfies the 
first $\ell$ constraints of \eref{pprime}, 
then $X \in 0 \oplus \psd{n-r}$ 
(this follows as in the proof of the ``if'' direction in Theorem \ref{ref-SDP}). 

Consider the homogenized SDP and its dual 
$$
\ba{rrrclcrrcll}
&\sup &  x_0 & &  &\clps \hspace{-0.4cm} &\inf & 0 & & &\\
\phom &s.t. & A_{i}^{\prime}\bullet X - b_{i}^{\prime} x_0 & = & 0 \,\, \forall i  &\clps\hspace{-0.4cm} & s.t. & \sum_{i}y_{i}A_{i}^{\prime} & \in & K^* & \dhom \\
 &  &- X &\leq_K & 0    & \clps \hspace{-0.4cm} &  &\sum_i y_{i}b_{i}^{\prime} & =    &-1.& 
\ena
$$
The optimal value of $\phom$ is $0, \,$ since if $(X, x_0)$ were feasible in it with $x_0 > 0, \,$ then 
$(1/x_0)X$ would be feasible in \eref{pprime}. 

We first check whether $\phom$ is strictly feasible, by solving the primal-dual pair of auxiliary SDPs
$$
\begin{array}{lrrllcrrcll}
    &   \sup  & t     &  &       & \clpss  &  \inf    &  0   &  & &  \\                   
\paux &   s.t.  & A_{i}^{\prime}\bullet X - b_{i}^{\prime} x_0 & = & 0 \,\, \forall i & \clpss &  s.t. &  \sum_i y_i A_i^{\prime} & \in & K^*   & \daux \\
    &   & -X + t \bpx 0 & 0 \\ 0 & I \epx  & \leq_K & 0 & \clpss &   &   \sum_i y_i b_i^{\prime} & = & 0 &   \\
    &         &     & & & \clpss &  & (\sum_i y_i A_i^{\prime}) \bullet \bpx 0 & 0 \\ 0 & I \epx & = & 1. 
\end{array}                                                                                              
$$

Clearly, $\paux$ is 
strictly feasible, with $(X, x_0, t) = (0, 0, -1)$ so it has strong duality 
with $\daux.$ Therefore 
\beqast
\phom \text{ is not strictly feasible} & \LRA & \text{the value of } \paux \, \text{ is }  0 \\
                                                    & \LRA & \paux \text{ is bounded} \\
						    & \LRA & \daux \text{ is feasible}. 
\eeqast
We distinguish two cases:

\nin {\bf Case 1:} 
$n-r \geq 2$ and $\phom$ is not strictly feasible. 

\begin{figure}[ht]
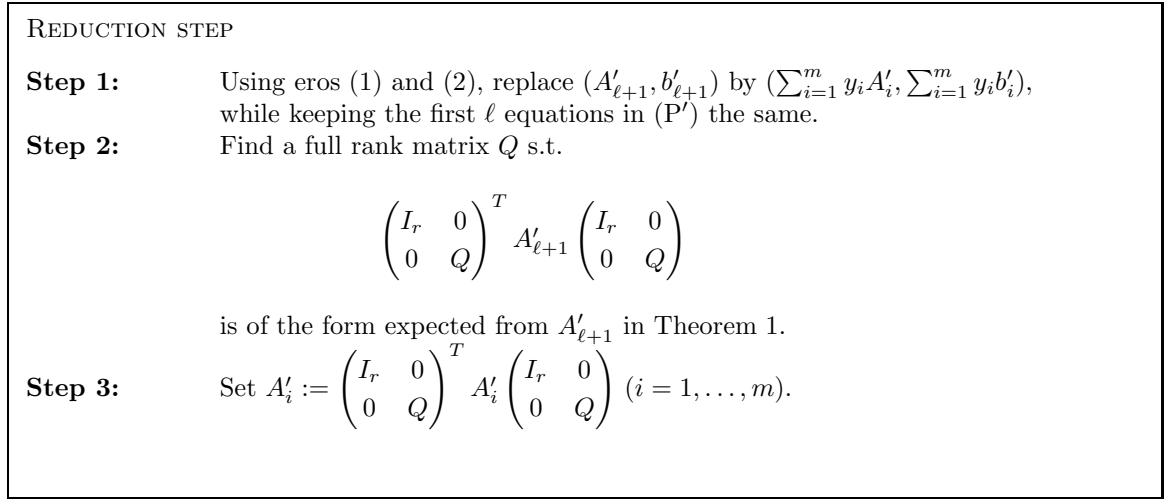

\framebox[6.05in]{\parbox{5.85in}{ 
{\sc Reduction step}
\begin{tabbing}{}
**************\=******\=*******\=***\= \hspace{2.5in} \=  \kill  
{\bf Step 1:} \> Using eros (1) and (2), replace $(A_{\ell+1}^\prime, b_{\ell+1}^\prime)$ by $(\sum_{i=1}^{m} y_{i}A_{i}^\prime, \sum_{i=1}^{m}y_{i}b_{i}^\prime), $ \\
               \> while keeping the first $\ell$ equations in \eref{pprime} the same. \\
{\bf Step 2:}  \> Find a full rank matrix $Q$ s.t. \\ \\
               \> \hspace{2cm} $\bpx I_r & 0 \\ 0 & Q \epx^T A_{\ell+1}^{\prime} \bpx I_r & 0 \\ 0 & Q \epx$ \\ \\
               \> is of the form expected from $A_{\ell+1}^\prime$ in Theorem \ref{ref-SDP}. \\
{\bf Step 3:}  \> Set $A_{i}^{\prime}:= \bpx I_r & 0\\ 0 & Q\epx ^{T} A_{i}^{\prime}\bpx I_r & 0\\ 0 & Q\epx\, (i=1, \dots, m).$ \\
\end{tabbing}
}}
\caption{The reduction step used to obtain the reformulated system} 
\label{fig:fr-alg}
\end{figure}

Let $y$ be a feasible solution of 
$\daux$ and apply the reduction step in Figure \ref{fig:fr-alg}. Now the lower 
$(n-r)$ by $(n-r)$ 
block of $\sum_i y_i A_i^{\prime}$ is nonzero, hence after Step 3 we have 
$r_{\ell+1}>0.$ We then set $\ell = \ell+1, \,$ and continue.

\vspace{.35cm}

\nin {\bf Case 2:} 
$n - r \leq 1$ or $\phom$ is strictly feasible. 

Now strong duality holds between $\phom$ and $\dhom$; when $n-r \leq 1, \,$ this is true because then 
$K$ is polyhedral. Hence $\dhom$ is feasible. Let $y$ be feasible in $\dhom$ and apply the same reduction step in Figure 
\ref{fig:fr-alg}. Then we set $k=\ell, \,$ and stop with the reformulation \eref{pref}.

We now complete the correctness proof of the algorithm. 
First, we note that the choice of the rotation matrix in Step 2 of the reduction steps implies 
that $A_1^\prime, \dots, A_\ell^\prime$ remain in the required form: cf. equation \eref{rotatez}. 

Second, we prove that after finitely many steps our algorithm ends in Case 2.
In each iteration both $\ell$ and $r = r_{1}+\ldots+r_{\ell}$ increase. 
If $n-r$ becomes less than or equal to $1$, then our claim is obviously true.
Otherwise, at some point during the algorithm we find $\ell=m-1.$ 
Then $b_{m}^{\prime}\neq 0, \,$ since \eref{pprime} is 
infeasible. 
Hence for any $X \in 0\oplus \pd{n-r}$ we can choose 
$x_0$ to satisfy the last equality constraint of $\phom$, hence at this point we are in Case 2.
\qed

We next illustrate our algorithm: 

\bex{\rm \label{ex1}
Consider the semidefinite system with $m = 6, \,$ and data
\begin{small}
\beqast
\ba{rclclc}
\vspace{0.2cm}
A_{1}=&\bpx 2&0&0&1\\0&3&0&-1\\0&0&4&2\\1&-1&2&0\epx, & A_{2}=&\bpx -1&2&1&-2\\2&3&3&1\\1&3&4&-3\\-2&1&-3&3\epx, & A_{3}=&\bpx -1&1&-2&0\\1&-2&0&2\\-2&0&-3&-2\\0&2&-2&-1\epx,\\
\vspace{0.2cm}
A_{4}=&\bpx 0&0&1&0\\0&1&0&0\\1&0&-1&0\\0&0&0&1\epx, & A_{5}=&\bpx 0&-1&0&0\\-1&0&0&-1\\0&0&1&1\\0&-1&1&0\epx, &A_{6}=&\bpx -1&0&0&-1\\0&0&1&0\\0&1&0&-1\\-1&0&-1&1\epx,\\
b=&(0,6,-3,2,1,3).
\ena
\eeqast
\end{small}

In the first iteration we are in Case 1, and find
\begin{small}
\beqast
y & = & (1,-1,-1,-1,-4,3),  \\
\sum_{i}y_{i}A_{i}^\prime &=& \bpx 1 & 1 & 0 & 0\\1&1&0&0\\0&0&0&0\\0&0&0&0\epx, \\
\sum_{i}y_{i}b_{i} &= & 0.  
\eeqast
\end{small}
We choose 
\begin{small}
$$
Q=\bpx 1&-1&0&0\\0&1&0&0\\0&0&1&0\\0&0&0&1 \epx
$$
\end{small}
to diagonalize $\sum_{i}y_{i}A_{i}^\prime, \,$ and after the reduction step 
we have a reformulation with data 
\begin{small}
\beqast
\ba{rclclc}
\vspace{0.2cm}
A_{1}^{\prime}=&\bpx 1&0&0&0\\0&0&0&0\\0&0&0&0\\0&0&0&0\epx, & A_{2}^{\prime}=&\bpx -1&3&1&-2\\3&-2&2&3\\1&2&4&-3\\-2&3&-3&3\epx, & A_{3}^{\prime}=&\bpx -1&2&-2&0\\2&-5&2&2\\-2&2&-3&-2\\0&2&-2&-1\epx,\\
\vspace{0.2cm}
A_{4}^{\prime}=&\bpx 0&0&1&0\\0&1&-1&0\\1&-1&-1&0\\0&0&0&1\epx, & A_{5}^{\prime}=&\bpx 0&-1&0&0\\-1&2&0&-1\\0&0&1&1\\0&-1&1&0\epx, &A_{6}^{\prime}=&\bpx -1&1&0&-1\\1&-1&1&1\\0&1&0&-1\\-1&1&-1&1\epx, \\
b^{\prime}=&(0,6,-3,2,1,3). 
\ena
\eeqast
\end{small}
We start the next 
iteration with this data, and $\ell = 1, \, r_1 = r=1. \,$ 
We are again in Case 1, and  find
\begin{small}
\beqast 
y & = & (0,1,1,0,3,-2), \\ 
\sum_i y_i A_i^\prime & = &\bpx 0 & 0 & -1 & 0\\0&1&2&0\\-1&2&4&0\\0&0&0&0\epx, \\
\sum_i y_i b_i^\prime & = & 0.
\eeqast
\end{small}
Now the lower right $3$ by $3$ block of $\sum_i y_i A_i^\prime$ is psd, and rank 1. 
We choose 
\begin{small}
$$
Q = \bpx 1&-2&0\\0&1&0\\0&0&1 \epx
$$
\end{small}
to diagonalize this block, 
and after the reduction step we have a reformulation with data 
\begin{small}
\beqast
\ba{rclclc}
\vspace{0.2cm}
A_{1}^{\prime}=&\bpx 1&0&0&0\\0&0&0&0\\0&0&0&0\\0&0&0&0\epx,& A_{2}^{\prime}=&\bpx 0&0&-1&0\\0&1&0&0\\-1&0&0&0\\0&0&0&0\epx, & A_{3}^{\prime}=&\bpx -1&2&-6&0\\2&-5&12&2\\-6&12&-31&-6\\0&2&-6&-1\epx,\\
\vspace{0.2cm}
A_{4}^{\prime}=&\bpx 0&0&1&0\\0&1&-3&0\\1&-3&7&0\\0&0&0&1\epx,& A_{5}^{\prime}=&\bpx 0&-1&2&0\\-1&2&-4&-1\\2&-4&9&3\\0&-1&3&0\epx, & A_{6}^{\prime}=&\bpx -1&1&-2&-1\\1&-1&3&1\\-2&3&-8&-3\\-1&1&-3&1\epx,\\
b^{\prime}=&(0,0,-3,2,1,3).
\ena
\eeqast
\end{small}
We start the last iteration with 
$ \ell = 2, \, r_1 = r_2 = 1, \, r = 2. \, $ 
We end up in Case 2, with 
\begin{small}
\beqast
y & = & (0,0,1,2,1,-1), \\
\sum_{i}y_{i}A_{i}^{\prime} &= &\bpx 0 & 0 & 0 & 1\\0&0&-1&0\\0&-1&0&0\\1&0&0&0\epx, \\
\sum_{i}y_{i} b_{i}^{\prime}&=&-1.
\eeqast
\end{small}
Now the lower right 2 by 2 submatrix of $\sum_{i}y_{i}A_{i}^{\prime}$ 
is zero, so we don't need  to rotate. After the reduction step 
the data of the final reformulation is
\begin{small}
\beqast
\ba{rclclc}
\vspace{0.2cm}
A_{1}^{\prime}=&\bpx 1&0&0&0\\0&0&0&0\\0&0&0&0\\0&0&0&0\epx, & A_{2}^{\prime}=&\bpx 0&0&-1&0\\0&1&0&0\\-1&0&0&0\\0&0&0&0\epx, & A_{3}^{\prime}=&\bpx 0&0&0&1\\0&0&-1&0\\0&-1&0&0\\1&0&0&0\epx,\\
\vspace{0.2cm}
A_{4}^{\prime}=&\bpx 0&0&1&0\\0&1&-3&0\\1&-3&7&0\\0&0&0&1\epx, &A_{5}^{\prime}=&\bpx 0&-1&2&0\\-1&2&-4&-1\\2&-4&9&3\\0&-1&3&0\epx, & A_{6}^{\prime}=&\bpx -1&1&-2&-1\\1&-1&3&1\\-2&3&-8&-3\\-1&1&-3&1\epx,\\
b^{\prime}=&(0,0,-1,2,1,3).
\ena
\eeqast
\end{small}
}
\eex

\section{The elementary reformulation of feasible systems} \label{feasible} 

For feasible systems we have the following result:

\bth \label{ref-SDP-2} Let $p \geq 0$ be an integer. Then the following hold: 
\benum
\item[(1)] 
The system \eref{p} is feasible with a maximum rank solution of rank $p$ 
if and only if it has a a feasible solution with rank $p$ and an elementary reformulation 
\begin{equation}\label{pref-feas}\tag{P$_{\rm ref,feas}$}
\begin{array}{rcl} 
 A_{i}^{\prime}\bullet  X & = & 0 \, (i=1,\ldots,k) \\
 A_{i}^{\prime}\bullet X & = & b_{i}^{\prime} \, (i=k+1,\ldots,m) \\
 X & \succeq & 0, 
\end{array}
\end{equation}
where $A_1^\prime, \dots, A_k^\prime$ are as in Theorem \ref{ref-SDP},
$$
k \geq 0, \, r_1 > 0, \dots, r_k > 0, \, r_1 + \dots + r_k =n-p, 
$$
and matrices $A_{k+1}^\prime, \dots, A_{m}^\prime$ and scalars 
$b_{k+1}^\prime, \dots, b_{m}^\prime$ are arbitrary. 
\item[(2)] 
Suppose that \eref{p} is feasible. Let \eref{pref-feas} be as above, and 
$({\rm P}_{\rm ref,feas, red})$ the system 
obtained from it by replacing the constraint 
$X \succeq 0$ by \mbox{$X \in 0 \oplus \psd{p}$}. 
Then
$({\rm P}_{\rm ref,feas, red})$ is well-behaved, i.e., for all $C \in \sym{n}$ the SDP
\beq \label{prefred}
\sup \, \{ \, C \bullet X \, | \, X \, {\rm \, is \, feasible \, in \,} ({{\rm P}}_{\rm ref,feas, red})\, \}
\eeq
has strong duality with its Lagrange dual 
\beq \label{drefred}
\inf \, \{ \, \sum_{i=1}^m y_i b_i \, : \, \sum_{i=1}^m y_i A_i^\prime - C \in (0 \oplus \psd{p})^* \, \}.
\eeq
\eenum
\enth

Before the proof we remark that 
the case $k=0$ corresponds to \eref{p} being strictly feasible.

\pf{of ``if'' in (1)} 
This implication follows similarly as in Theorem \ref{ref-SDP}.

\pf{of (2)} 
This implication follows, since 
$({{\rm P}}_{\rm ref,feas, red})$ is trivially strictly feasible.

\pf{of ``only if'' in (1)}  
We modify the algorithm that we used to 
prove Theorem \ref{ref-SDP}. 
We now do not assume that \eref{p} is infeasible, nor that the optimal value of $\phom$ is zero.
As before, we keep iterating in Case 1, until we end up in Case 2, with strong 
duality between $\phom$ and $\dhom.$ We distinguish two subcases: 

\nin {\bf Case 2(a):} 
The optimal value of $\phom$ is $0. \,$ We proceed as before to construct the $(k+1)^{st}$ equation 
in \eref{pref}, which proves infeasibility of \eref{pprime}. 

\nin {\bf Case 2(b):} 
The optimal value of $\phom$ is positive (i.e., it is $+ \infty$). 
We choose 
$$
(X, x_0) \in K \times \rad{}
$$
to be feasible, with $x_0 > 0.$ Then $(1/x_0)X$ is feasible in \eref{pprime}, but it may not have 
maximum rank.
We now construct a maximum rank feasible solution in \eref{pprime}. If 
$n-r \leq 1, \,$ then a simple case checking can complete the construction.
If $n-r \geq 2, \,$ then we take 
$$
(X', x_0') \, \in \, \ri K \times \rad{} 
$$
as a strictly feasible solution of $\phom.$ Then for a small $\epsilon > 0$ we have that 
$$
(X + \eps X', x_0 + \eps x_0') \, \in \, \ri K \times \rad{}
$$
is feasible in $\phom$ with $x_0 + \eps x_0'  > 0.$ Hence 
$$
\dfrac{1}{x_0 + \epsilon x_0'} (X + \eps X') \in \ri K
$$
is feasible in \eref{pprime}.
\qed

\bex{\rm \label{ex2}
Consider the feasible semidefinite system with $m = 4, \,$ and data
\beqast
\ba{rlcllc}
\vspace{0.2cm}
A_{1}=&\bpxastr -2&2&7&-3\\2&-2&-4&-6\\7&-4&-15&-7\\-3&-6&-7&0\epxastr,& A_{2}=&\bpxastr 2&0&-3&2\\0&4&6&4\\-3&6&14&5\\2&4&5&0\epxastr,& A_{3}=&\bpxastr 2&0&-3&-1\\0&-1&-3&0\\-3&-3&-3&2\\-1&0&2&0\epxastr,\\
\vspace{0.2cm}
A_{4}=&\bpxastr -1&1&4&2\\1&6&11&2\\4&11&16&1\\2&2&1&0\epxastr ,& b = &(-3,2,1,0). \,& 
\ena
\eeqast
The conversion algorithm produces the following $y$ vectors, and rotation matrices: it produces 
\beq
\ba{rlllllll} 
y & = & (1, 2, -1, -1), & V & = & \bpxastr 1 & -1 & 0 & 0 \\ 0 & 1 & 0 & 0 \\ 0 & 0 & 1 & 0 \\ 0 & 0 & 0 & 1 \epxastr \\
\ena
\eeq
in step 1, and 
\beq
\ba{rllllllllllllll}
y & = & (0, 1, -2, -1), & V & = & \bpxastr 1 & 0 & 0 & 0 \\ 0 & 1 & -2 & 0 \\ 0 & 0 & 1 & 0 \\ 0 & 0 & 0 & 1 \epxastr
\ena
\eeq
in step 2 (for brevity, we now do not show the $\sum_i y_i A_i$ matrices, and the intermediate data). 
We obtain an elementary reformulation with data and maximum rank feasible solution 
\beqast
\ba{rlllll}
\vspace{0.2cm}
A_{1}^{\prime}=&\bpxastr 1&0&0&0\\0&0&0&0\\0&0&0&0\\0&0&0&0\epxastr ,& A_{2}^{\prime}=&\bpxastr -1&0&-1&2\\0&1&0&0\\-1&0&0&0\\2&0&0&0\epxastr,& A_{3}^{\prime}=&\bpxastr 2&-2&1&-1\\-2&1&-2&1\\1&-2&1&0\\-1&1&0&0\epxastr,\\
\vspace{0.2cm}
A_{4}^{\prime}=&\bpxastr -1&2&0&2\\2&3&1&0\\0&1&0&1\\2&0&1&0\epxastr ,& b^{\prime}=&(0,\;0,\;1,\;0), & X=& \bpxastr 0&0&0&0\\0&0&0&0\\0&0&1&0\\0&0&0&1 \epxastr.
\ena
\eeqast
}
\eex

In the final system the first two constraints prove that the rank of 
any feasible solution is at most $2.$ Thus the system itself and $X$ are a 
certificate that $X$ has maximum rank, hence 
it is easy to convince a ``user'' that $({\rm P}_{\rm ref,feas, red})$ (with $p=2$) is 
strictly feasible, hence well behaved. 

\section{Discussion} \label{conclusion}

In this section we discuss our results in some more detail.

We first compare our conversion algorithm with facial reduction algorithms,
and describe how to adapt the algorithm of Waki and Muramatsu \cite{WakiMura:12} to obtain our 
reformulations. 
\begin{Remark}{\rm 
We say that a convex subset $F$ of a convex set $C$ is a {\em face of $C$} , if
$x, y \in C, \, 1/2(x+y) \in F$ implies that $x$ and $y$ are in $F.$
When \eref{p} is feasible, 
we define its {\em minimal cone} as the smallest face of $\psd{n}$ that contains 
the feasible set of \eref{p}. 

The algorithm of Borwein and Wolkowicz \cite{BorWolk:81, BorWolk:81B} 
finds the minimal cone of a feasible, but possibly nonlinear 
conic system. The algorithm of Waki and Muramatsu \cite{WakiMura:12} is a 
simplified variant which is applicable to conic linear systems, and can detect infeasibility.
We now describe their Algorithm 5.1, which specializes their general algorithm 
to SDPs, and how to modify it to obtain our reformulations. 

In the first step they find $y \in \rad{m}$ with 
$$
W \, := \, - \sum_{i=1}^m y_i A_i \succeq 0, \, \sum_{i=1}^m y_i b_i \geq 0.
$$
If the only such $y$ is $y=0, \,$ they stop with $F = \psd{n}; $ 
if $\sum_{i=1}^m y_i b_i > 0, \, $ they stop and report that \eref{p} is infeasible.
Otherwise they replace $\psd{n}$ by 
$\psd{n} \cap W^\perp,$ apply a rotation step 
to reduce the order of the SDP to $n-r, \, $ where $r$ is the rank of $W,$ 
and continue. 

Waki and Muramatsu do not apply elementary row operations. We can obtain our reformulations 
from their algorithm, if after each iteration $\ell = 0, 1, \dots$ we
\bit
\item choose the rotation matrix to turn the psd part of $W$ into $I_{r_\ell}$ for some 
$r_\ell \geq 0.$ 
\item add eros to produce an equality constraint like 
the $\ell$th constraint in \eref{pref}, or \eref{pref-feas}.
\eit
In their reduction step they also rely on Theorem 20.2 from Rockafellar 
\cite{Rockafellar:70}, while we use explicit SDP pairs. 
For an alternative approach to ensuring strong duality, called {\em conic expansion}, 
we refer to Luo et al \cite{LuoSturmZhang:97}; 
and to \cite{WakiMura:12} for a detailed study of the connection of the two approaches. 
}
\end{Remark}

We next comment on how to find the optimal solution of a linear function over the 
original system \eref{p}, and on duality properties of this system.
\brem{\rm 
Assume that \eref{p} is feasible, and we used 
the rotation matrix $V$ to obtain 
\eref{pref-feas} from \eref{p}. Let $C \in \sym{n}.$ 
Then one easily verifies 
\beqast
\sup \, \{ \, C \bullet X \, | \, X \, {\rm \, is \, feasible \, in \,} \eref{p}\,\} & = & 
\sup \, \{ \, V^{T} C V \bullet X \, | \, X \, {\rm \, is \, feasible \, in \,} (\rm {P}_{\rm ref,feas}) \, \} \\
  & = & \sup \, \{ \, V^{T} C V \bullet X \, | \, X \, {\rm \, is \, feasible \, in \,} ({{\rm P}}_{\rm ref,feas, red})\, \},
\eeqast
and by Theorem \ref{ref-SDP-2} the last SDP has strong duality with its Lagrange dual. 

Clearly, \eref{p} is well behaved, if and only if its ESD-reformulations are. 
The system \eref{p}, or equivalently, system \eref{pref-feas} 
may not be well behaved, of course. 
We refer to \cite{Pataki:10} for an exact characterization of well-behaved semidefinite systems
(in an inequality constrained form). 
}
\erem 

We next comment on algorithms to generate the data of all SDPs which are either infeasible,
or have a maximum rank solution with a prescribed rank. 
\brem{\rm Let us fix an integer $p \geq 0, \,$  and define the sets 
%
\beqast
\INFEAS  & = & \{ \, (A_i, b_i)_{i=1}^m \, \in (\sym{n} \times \rad{})^m \, : \, \eref{p} \, \mathrm{is \; infeasible \,} \}, \\
\FEAS(p) & = & \{ \, (A_i, b_i)_{i=1}^m \,  \in (\sym{n} \times \rad{})^m \, : \, \eref{p} \, \mathrm{is \; feasible,} \mathrm{\;with \; maximum}\\
& & \mathrm{\hspace{0.3cm} rank\; solution \; of \; rank \; } p \, \}.
\eeqast

These sets -- in general --  are nonconvex, neither open, nor closed.
Despite this, we can systematically generate {\em all} of their elements.
To generate all elements of $\INFEAS,$ we use Theorem \ref{ref-SDP}, by which 
we only need to find  systems
of the form \eref{pref}, then reformulate them. 
To generate all elements of $\FEAS(p)$ we first find 
constraint matrices in a system like \eref{pref-feas}, then choose 
$X \in 0 \oplus \pd{p}, \,$ and set $b_i^\prime := A_i^\prime \bullet X \, $ for all $i.$ 
By Theorem \ref{ref-SDP-2} all elements of $\FEAS(p)$ 
arise as a reformulation of such a system. 

Loosely speaking, Theorems \ref{ref-SDP} and \ref{ref-SDP-2} show that there are only finitely many ``schemes'' to 
generate an infeasible semidefinite system, and  a feasible system with a maximum rank solution having a prescribed rank. 

The paper \cite{Pataki:10} describes a systematic method to generate all well behaved semidefinite systems 
(in an inequality constrained form), in particular, to generate all linear maps under which 
the image of $\psd{n}$ is closed. Thus, our algorithms to generate 
$\INFEAS$ and $\FEAS(p)$ complement the results of \cite{Pataki:10}. 
}
\erem

We next comment on strong infeasibility of \eref{p}. 
\brem{\rm 
Clearly, \eref{p} is strongly infeasible (i.e., \eref{palt} is feasible), if and only if 
it has a reformulation of the form \eref{pref} with $k = 0.$ 
Thus we can easily generate the data of all strongly infeasible SDPs: we only need to find
systems of the form \eref{pref} with $k=0, \,$ then reformulate them. 

We can also easily generate weakly infeasible instances using Theorem \ref{ref-SDP}: we can choose
$k+1 = m, \,$ and suitable blocks of the $A_i^\prime$ in \eref{pref} to make sure that they do not have a psd linear 
combination. (For instance, choosing the block of $A_{k+1}^\prime$ that corresponds to rows 
$r_1 + \dots + r_{k-1} +1$ through $r_1 + \dots + r_k$ and the last $n - r_1 - \dots - r_{k+1}$ columns will do.)
Then \eref{pref} is weakly infeasible. It is also likely to be weakly infeasible, if we 
choose the $A_i^\prime$ as above, and $m$ only slightly larger than $k+1.$ 

Even when \eref{p} is strongly infeasible, our conversion algorithm may only find a reformulation with $k>0. $ 
To illustrate this point, consider the system with data 
\beq
\ba{rclclc}
\vspace{0.2cm}
A_{1}=&\bpx 1&0&0\\0&0&0\\0&0&0\epx, & A_{2}=&\bpx 0&1&0\\1&0&0\\0&0&0\epx, & A_{3}=&\bpx 0&0&0\\0&1&0\\0&0&1\epx, \\
b =& (0, -1, 1).
\ena
\eeq
This system is strongly infeasible (\eref{palt} is feasible with $y = (4,2,1)$), and it is already in the 
form of \eref{pref} with $k=1. \,$ Our conversion algorithm, however, 
constructs a reformulation with $k=2,$ 
since  it finds $\phom$ to be not strictly feasible in the first two steps. 
}
\erem

We next discuss complexity implications. 
\brem{\rm 
Theorem \ref{ref-SDP} implies that semidefinite feasibility is in 
$\conp$ in the real number model of computing.
This result was already proved by Ramana \cite{Ramana:97} and 
Klep and Schweighofer \cite{KlepSchw:12} via 
their Farkas' lemma that relies on extra variables. To check the 
infeasibility of \eref{p} using our methods, 
we need to verify that \eref{pref} is a 
reformulation of \eref{p}, using eros, and a rotation matrix $V.$ 
Alternatively, one can check that 
$$
A_i^\prime \, = \, V^T \bigl( \sum_{j=1}^m t_{ij} A_j)V \, (i=1, \dots, m) 
$$
holds for $V$ and an invertible matrix $T = (t_{ij})_{i,j=1}^m.$ 
}
\erem

\section{Conclusion}

Two  well-known pathological phenomena in semidefinite programming are 
that Farkas' lemma may fail to prove infeasibility, and strong duality does not hold in general. 
Here we described an exact certificate of infeasibility, and a strong dual for SDPs, which do not assume 
any constraint qualification. 
Such certificates and duals have been known before: see 
\cite{BorWolk:81, BorWolk:81B, Ramana:97, RaTuWo:97, WakiMura:12, KlepSchw:12, Pataki:13}. 

Our approach appears to be simpler: in particular, the validity of our infeasibility certificate 
-- the infeasibility of the system \eref{pref} -- is almost a tautology 
(we borrow this terminology from the paper \cite{Papp:Ali} on semidefinite 
representations). We can also easily 
convince a ``user'' that the system $({\rm P}_{\rm ref,feas, red})$ is well behaved
(i.e., strong duality holds for all objective functions).
To do so, we use a maximum rank feasible solution, and the system itself, which 
proves that  this solution has maximum rank. 

In a somewhat different language, elementary reformulations provide a standard form of spectrahedra -- the feasible sets 
of SDPs -- to easily check their emptiness, or a tight upper bound on the rank of feasible solutions.
We hope that these standard forms will be useful in studying the geometry of spectrahedra -- a subject of 
intensive recent research \cite{NetzerPlaumSchwei:10, Blekhetal:12, Vinzant:14, SinnSturmfels:14}. 

\nin{\bf Acknowledgement} We thank Rekha Thomas and the anonymous referees for their careful reading of the paper,
and their constructive comments.


\bibliography{mysdp}

\end{document}